 \documentclass[article,12pt]{elsarticle}
 
\usepackage{graphicx} 
 
\usepackage[colorlinks=true, allcolors=blue]{hyperref}

\usepackage{amssymb}
\usepackage{amsmath}
\usepackage{amsthm}

\newcommand{\cS}{\cal S}

 \usepackage{bm}

\newtheorem{definition}{Definition}

\newtheorem{remark}{Remark}

\newtheorem{assumption}{Assumption}
\newtheorem*{hypotheses*}{Hypotheses}

\usepackage{algorithm}%
\usepackage{algorithmicx}%
\usepackage{algpseudocode}%
\usepackage{listings}%

 \usepackage{verbatim}
\usepackage{epsfig}

 \usepackage{caption}
\usepackage{xcolor}
\usepackage{bm}
\usepackage{multirow}

\usepackage[normalem]{ulem}
\usepackage{lineno}
\usepackage{ subfigure }

\journal{Journal of Computational and Applied Mathematics}

\begin{document}

\begin{frontmatter}



\title{A regularization method  for planar offset curves and    {bi-offset} recognition}


\author[label1]{Rosanna Campagna}
\author[label2]{Salvatore Mondrone}
\author[label3]{Tomas Sauer} 

\affiliation[label1]{organization={Department of Mathematics and Physics, University of Campania ``L. Vanvitelli''},
            state={Italy}}
\affiliation[label2]{organization={Department of Information Engineering, Computer Science and Mathematics, University of L'Aquila},
            state={Italy}}
\affiliation[label3]{organization={Lehrstuhl f\"ur Mathematik
	mit Schwerpunkt Digitale Bildverarbeitung \& FORWISS, Universit\"at
	Passau},
            state={Germany}}

\begin{abstract}
Offset curves for planar trajectories are interesting in the generation of tool paths for numerically controlled industrial machines and  in trajectory planning methods for autonomous driving systems. 
Theoretical offset curves  
may exhibit peculiar singularities, including self-intersections, which limit their use in practical applications. Existing approaches   address these issue through geometric filtering techniques to detect and remove undesirable features  but the computation of accurate and well-behaved offset curves remains a challenging task.
We assume a first stage of functional approximation of trajectories by penalized Hermite spline regression  enabling the simultaneous fitting of positions and tangents. The regularization is imposed on the second derivatives, effectively mitigating the jerk effect, which is particularly relevant in motion planning and path smoothing applications. 
 Then, taking into account the geometrical pointwise properties of the resulting curve, we design two offset curves through the simultaneous approximation of function values and derivatives.  
Then, a mathematical
model  to obtain the so-called \emph{bi-}offset as most fitting
as with the original generator curve is proposed, also relating 
  the offset range   and   pointwise curvature   values.  The adaptive reconstruction  of the center line from the external boundaries is a topic of interest 
  and is the main focus of our work. Numerical experiments confirm the reliability of our approach at every stage of the resolution process.
\end{abstract}




\end{frontmatter}




\section{Introduction}
 
\begin{figure}
  \centering \includegraphics[width=10cm]{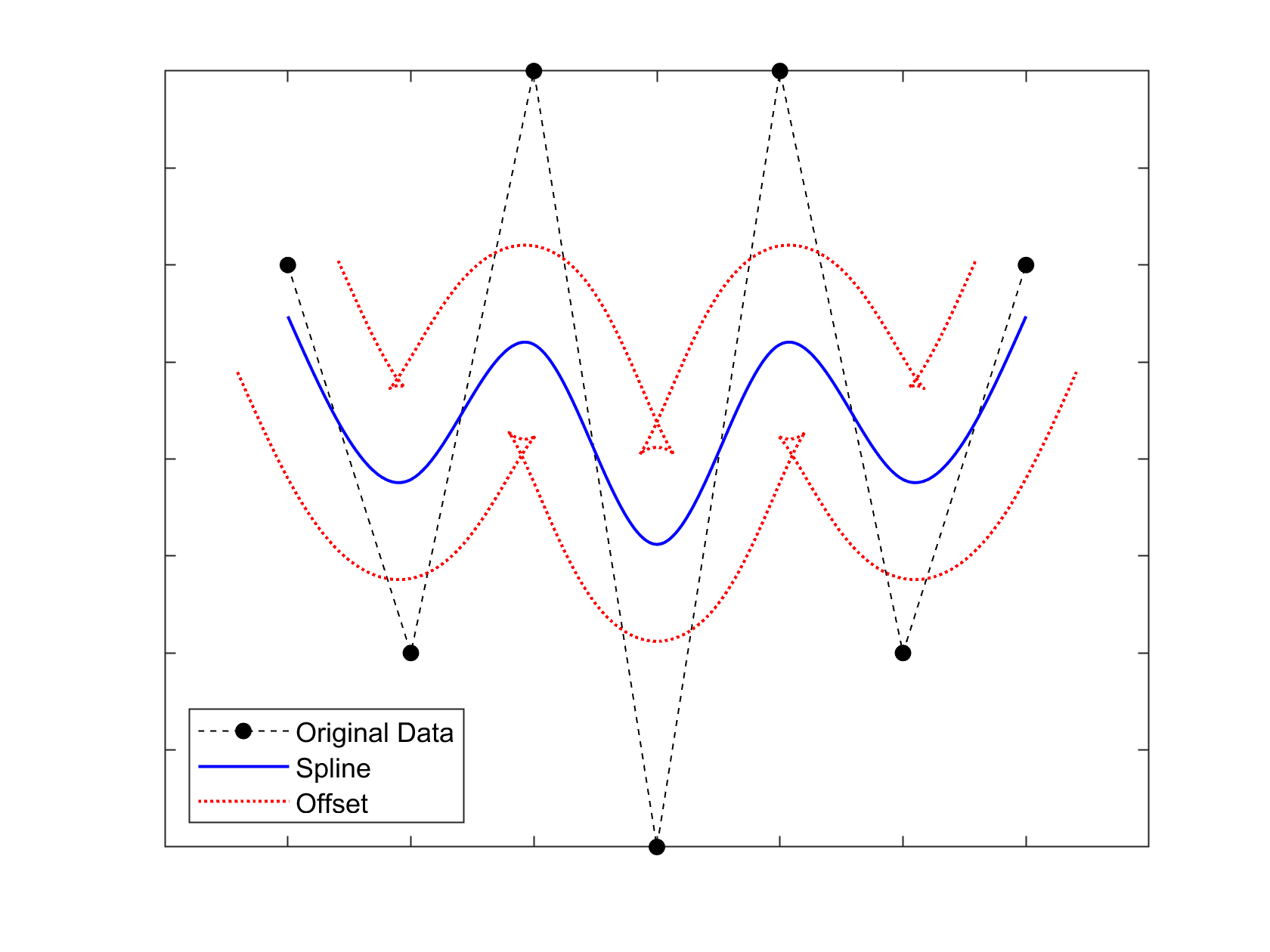}
  \caption{Cubic smoothing spline constructed on seven data points and its offsets}
\label{fig:offsetissue}
\end{figure}

The construction of offset curves for planar trajectories is a topic of widespread interest in numerous applications that require the determination of path approximations, the definition of edge trajectories, or the recognition and prediction of a path, starting from assigned or sensor-acquired control points. These examples include the generation of tool paths for
industrial numerical control machines and the development of trajectory planning
methods for autonomous driving systems. While the approximation
of planar curves using polynomial splines of arbitrary degree under geometric continuity
constraints has been widely studied, and several methodologies
for offset modeling are already available in the literature
 \cite{Hoschek1,HOSCHEK198759,Hoscheck2,HOSCHEK198827}, the computation of accurate and well-behaved offset curves and
 related recognition of paths   from their two offsets   remain   challenging tasks.
Offset curves typically do not belong to the same functional space as the original trajectory. For example, the offset of a spline curve is not generally a spline itself. Moreover, offset curves may exhibit peculiar singularities, including self-intersections, which complicate their use in practical applications (see Fig. \ref{fig:offsetissue}). 

Existing approaches in the literature address these issue through geometric filtering techniques to detect and remove undesirable features \cite{mca22030039}.
Particularly \cite{Hoschek1,Hoscheck2} address the offset problem,   accurately reproducing both the generating curve and its offsets through approximating curves belonging to the same functional class. However, they lack the \textit{reconstruction modeling}, which will be examined and developed in the present work.
 More into details, 
 {in~\cite{Hoschek1}, a sequence of points converging to the self-intersection of the offset curve is generated using Newton-like methods. Once the singular segment is identified in this manner, it is directly removed, thereby eliminating cusps and self-intersections.}
 {In~\cite{Hoscheck2}, the generator is a B\'ezier curve. By exploiting the geometric properties of spline curves and their geometric invariance, the offsets are approximated using spline curves by setting up a least-squares problem. This formulation is derived from imposing geometric continuity conditions. The selection of the optimal parameters for the approximation follows the approach of Hoschek's studies cited in~\cite{HOSCHEK198759,HOSCHEK198827}.
}
 {In~\cite{FAROUKI199083}, the main analytical-differential properties of offset curves are outlined. In particular, characterizations of cusp points in relation to the evolute of the original curve are provided, along with the regularization property of trimmed offsets in terms of the energy integral, and the mutual distance between the offset and the generator curve.}
 {In~\cite{FAROUKI1990101}, a purely algebraic study of offset curves is conducted for the case of polynomial or rational generator curves. Recalling the notion of the resultant of polynomials, it is applied to construct the implicit equations of the offset curves and to study their respective degree. A similar analysis is used to identify self-intersection points, characterized as zeros of the determinant of a modified Sylvester matrix.}

In \cite{Joon}, a Hermite interpolation scheme is proposed to approximate offset curves when the original curve is   regular, producing a polynomial representation of the offset. However, such direct approximations often fail to preserve essential characteristics of the original curve, due to the inherent nonlinearity of the offsetting operation, which involves a displacement along the normal vector field of the curve.


Recognizing these singular phenomena, which are critical in application contexts, focusing on  the analytical and algebraic properties of offset curves, we introduce     a new  numerical model and a related algorithm for designing offset curves that are regular, smoothing, \emph{well-performing}  and
useful to "coming back" to the original path.


Actually, drawing an    approximation  to  the  "offset of an offset"  curve,     is 
generally  no  longer  able  to  reproduce  the  center-line curve  from its boundaries, and  it  can  lead  to
very  unwanted  effects and artifacts.  The  reason  lies  in  the  nature  of  the
offset  formation,  in  which  the  shift  is  made  in  the  direction  of  the  normal  vector. 
However, if the derivatives of a curve are not precisely preserved,  no matter how well the values of the
curve may be approximated, the offset will be faulty.

In  this work, we present a regularization technique  {based on Hermite spline regression \cite{CCF}}, which incorporates both function values and derivative information, enabling the simultaneous fitting of positions and tangents. The regularization is imposed on the second derivatives, effectively mitigating the jerk effect, which is particularly relevant in motion planning and path smoothing applications. 
Taking into account the geometrical point-wise properties of the resulting curve, we design two offset curves through the simultaneous approximation of function values and derivatives; 
furthermore, the contribution of the first derivative in the fitting process is modulated by a weight parameter, which is heuristically determined to balance adherence to tangent data against overall offset curve smoothness. Finally,  
our method allows to reconstruct an approximation of the   base trajectory, given its offset curves (we say the bi-offsets, for short) that to the best of our knowledge  deserve   interest  and not so much    investigated in literature.

The paper is organized as follows:  in Section \ref{sec:2}, we introduce the definition of planar offset curve, giving basic    geometrical and algebraic properties; we also analyse the nature of some critical points of the offset and their strict correspondence with the analytical character of   a regular curve; 
Section \ref{sec:3} delineates a new algorithm that allows for the construction of offsets free from singularities. The proposed approach is configured with three main objectives: {removal of cusp and self-intersection points relying on a non-linear regression model; obtaining offset curves that may belong to the same functional class as the generator curve (an element that is generally lost); and finally, defining a mathematical model that allows to obtain  \textit{the  offset of the offset} (bi-offset) as most fitting as  with the original generator curve.} Some algorithmic details and remarks can be found in Section \ref{sec:4}. Numerical experiments are in Section \ref{sec:5}, giving results and comparative studies confirming the reliability of our approach.
Conclusions and future investigations close the work.

\section{Preliminaries\label{sec:2}}

To deal with  the offset definition and approximation, 
 it is necessary to recall some preliminaries.
 \begin{definition}
		Let $r:[0,1]\to \mathbb{R}^2$ be a differentiable curve in $\mathbb{R}^2$ defined on the compact interval $[0,1]$. Then $r$ is called a regular curve if
		\begin{equation*}
			\Vert r'(t) \Vert = \sqrt{x'^2(t) + y'^2(t)} \neq 0 \quad \text{for each}  \, \, t \in(0,1).
		\end{equation*}
\end{definition}
\begin{definition}\label{def:offsets}
Let $r:[0,1]\to \mathbb{R}^2$ be a regular 
curve, 
parametrically given  as 
$$r(t)=(x(t), \, y(t)), \quad t\in [a,b]$$
in the plane, and let $\tau \neq 0$ be a fixed signed distance. The {offset curve} generated by $r$ and located at a distance $\tau$ from it is defined as:
\begin{equation*}
	r_o(t)=r(t)+\tau N(t) \quad \text{for all} \,\, t\in[0,1], 
\end{equation*}
where
$$
N(t) =\frac{1}{\sqrt{x'(t)^2 + y'(t)^2}} (-y'(t), x'(t))
$$
is the principal unit normal field to $r(t)$ at $ t$.
 We obtain an exterior or interior offset according to whether $\tau>0$ or $\tau<0$.
\end{definition}
To every regular unit-speed planar curve, i.e., a curve with $\| r' \| = 1$, two unit vector fields are associated (namely the tangent and the normal), defined as follows:
\begin{equation}
T(t)=r'(t), \quad N(t)=\mathcal{J}r'(t)
\end{equation}
where the linear map $\mathcal{J}$, called the complex structure of $\mathbb{R}^2$, is
\begin{equation*}
\mathcal{J}: (p_1,p_2) \in \mathbb{R}^2 \mapsto (-p_2,p_1) \in \mathbb{R}^2.
\end{equation*}
\begin{definition}\label{curvature}
Let $r:[0,1]\to \mathbb{R}^2$ be a regular curve. The curvature $k$ of $r$ is given by the formula
\begin{equation}
	\label{fre or}
	k=\dfrac{r'' \cdot \mathcal{J}r'}{\Vert r' \Vert^3}
\end{equation}
\end{definition}
For any regular unit-speed curve in the plane, the Frenet equations are expressed as follows:
\begin{equation*}
T'=kN, \quad N'=-kT.
\end{equation*}

For our purpose, it is useful to compute the tangent and normal unit vector fields associated with the offset curve. We can express the successive derivatives of $r_o(t)$ as:
\begin{equation}
\begin{aligned}
	r_o'&=r'+\tau N'=r'+k\tau T=r'+k\tau r'=(1+k\tau )r' \\
	r_o''&=(1+k\tau)r''+\tau k'r',
\end{aligned}
\label{succ der off}
\end{equation}
where the quantities on the right-hand sides all pertain to the generator curve $r(t)$.
Hence, we can express the tangent and normal unit vector fields of the offset as:
\begin{equation}
T_o=\dfrac{r'_o}{\Vert r_o'\Vert}=\dfrac{(1+k\tau)}{\vert 1+k\tau \vert}T, \quad
N_o = \mathcal{J}\dfrac{r'_o}{\Vert r'_o\Vert} = \dfrac{(1+k\tau)}{\vert 1+k\tau \vert}N
\label{tan-nor-off}
\end{equation}
at the corresponding points of the offset and generator curves.
\begin{definition}
\label{def:ray-curvature}
Let $r:[0,1]\to \mathbb{R}^2$ be a regular unit-speed curve in the plane, let $k$ be its curvature, and let $r_o$ be the offset curve generated by $r$ at a fixed signed distance $\tau$. Then the offset curve $r_o$ exhibits an {ordinary cusp} at each parameter value $t_c$ such that:
\begin{equation*}
	k(t_c)=-\dfrac{1}{\tau} \quad \text{and} \quad k'(t_c) \neq 0.
\end{equation*}
If $k(t) \neq -\dfrac{1}{\tau}$ for all $t \in (0,1)$, the offset $r_o(t)$ at distance $\tau$ is said to be non-degenerate.
\end{definition}

We observe that, if $k(t_c)=-\dfrac{1}{\tau}$ and $k'(t_c) \neq 0$, then the factor $\dfrac{(1+\tau k)}{\vert 1+ \tau k \vert}$ in the expression of the tangent and normal field of the offset is a "step function" which changes abruptly from $-1$ to $+1$, or vice-versa, at $t_c$. The tangent $T_o$ and the normal $N_o$ suffer sudden inversions as we traverse such a point.

If $r_o$ has an ordinary cusp at a parameter value $t_c$, we take the tangent line to the offset at that cusp to be the line through $r_o(t_c)$ which is parallel to the tangent $T$ of the generator at the corresponding point $t_c$.

Now, substituting   (\ref{succ der off}) into the expression for curvature in (\ref{fre or}), we obtain:
\begin{equation}
\dfrac{ r''_o \cdot \mathcal{J}r'_o}{\Vert r'_o \Vert^3} =\dfrac{[(1+k\tau)r''+\tau k'r']\cdot(1+k\tau)\mathcal{J}r'}{(1+\tau k)^2\vert 1+\tau k\vert}= \dfrac{ k }{\vert 1+\tau k\vert}
\label{curv-off}
\end{equation}

\begin{remark}
When $k(t_e)=-\dfrac{1}{\tau}$ and $k'(t_e)$ vanishes but $k''(t_e)$ does not, the quantity
\begin{equation*}
	\dfrac{1+\tau k}{\vert 1+\tau k \vert}
\end{equation*}
has the same value on either side of $t_e$, and therefore $T_o$ and $N_o$ are not inverted on passing through $t_e$. Although $r_o$ has a continuous tangent at $t_e$, the expression for $k_o$ indicates that the offset curvature increases without bound.
\end{remark}

In this article, we   consider only a specific class of planar curves having the form $y=f(x)$. The function $f$ is assumed to be at least of class $\mathcal{C}^2$ on its domain. We recall that a Cartesian curve can be expressed equivalently as:
\begin{equation*}
r(t)=(t,f(t)) \quad \text{for all} \,\, t \in \text{dom}\, f.
\end{equation*}
Thus the offset curve and the curvature for a Cartesian curve are:
\begin{equation*}
r_o(t)=\left(t-\tau \dfrac{f'(t)}{\sqrt{1+f'(t)^2}},f(t)+\tau \dfrac{1}{\sqrt{1+f'(t)^2}} \right)
\end{equation*}
and
\begin{equation*}
k=\dfrac{f''(t)}{(1+f'(t)^2)^{3/2}}.
\end{equation*}

The next section presents the mathematical modeling of offset pairs starting from a functional fit of a given data set, correlated to the expected or predicted trajectory. The proposed approximation algorithm excludes the explicit computation of singular points and cusps. Instead, the procedure is framed within successive optimization problems, whose smooth regression solution is modeled through an appropriate constraint-driven minimization of the vertical distances between the points, enabling the adaptive construction of a  singularity-free curve.

\section{Offsets computation and trajectory recognition\label{sec:3}}

\begin{figure}
\centering \includegraphics[scale=0.2]{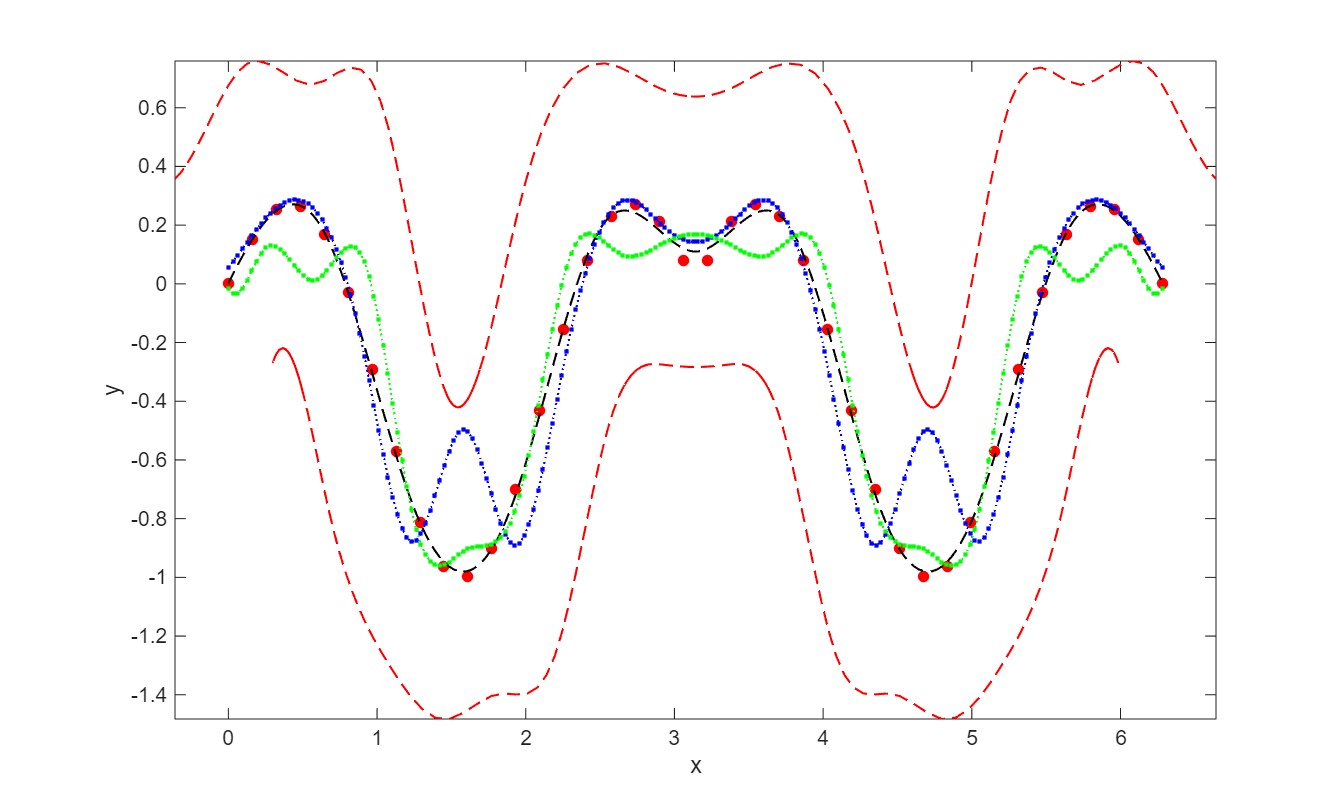}
\caption{Base  curve  (black '$\textcolor{red}{\bullet} -$'),    computed  offset   at distance $\tau=0.5$  (red '- -'), {interior bi-offset}   (blue '$\cdot$ -'),  {exterior bi-offset}   (green '$\cdot$ -').\label{fig:offsetpb}}  
\end{figure}

 
In this section we define an algorithm that, taking into account the limits of the offsets definition, gives a smooth and meaningful  offset curve, and also  allows to  "invert" this procedure, giving a functional form for the offset of the offset (we   say \emph{bi-offset}), as close as possible to the  original functional data-driven model.
One of the main issues  in the reconstruction of the offset of a planar curve is that the \textit{offset of an offset} results in a new curve, possibly far from the original generator.
 Fig. \ref{fig:offsetpb}  is to show  that reproducing the bi-offset, from above (\textit{interior bi-offset}) and from below (\textit{exterior bi-offset}), is generally not good for designing  the  base  curve.
 This issue arises from the definition of the offsets, 
where {the shift occurs along the normal vector,
  \emph{without preserving  the derivatives}}; so,  no matter how well the values of the
 curve may be approximated, {the offset will be faulty}.

 Based on these  considerations, we propose an algorithm designed to resolve all the aforementioned problems. The objective is to construct a model in which the offsetting process is \emph{reversible}.
 In order to conduct a pointwise analysis of the proposed algorithm, it is appropriate to decompose the overall procedure into a sequence of well-defined steps. At each step, a regression problem must be formulated, preceded by the construction of an appropriate minimization problem.  
 The numerical scheme will be synthesized in four steps; for any step further 
 algorithmic and computational  details   will be given in the next Section \ref{sec:4}.

 {To avoid redundancies in the next definitions, we preliminarily introduce some settings in the following Assumption.}
\begin{assumption}
\label{assump1}
Let consider a set of  $m$  data 
${\mathcal D}=\{(x_j,y_j)\}$, $j=1,\dots,m \in[a, b] \times \mathbb{R}\subseteq \mathbb{R}^{2} $  with $a=x_{1}<\cdots < b=x_{m}$ 
and a basis of $n$ B-splines $\{B_j\}_{j=1}^{n}$, of order $4$,
   defined on an augmented  set of $n+4$  uniformly distributed knots, 
  $\Xi=\{ \xi_{-1},\ldots,\xi_{n+2}\}$  with $h = \frac{b-a}{n-3}$, $\xi_{2} \equiv a$, $\xi_{n-1}\equiv b$, $\xi_{-1}< \xi_0 < \xi_1 < a$ and $b < \xi_n < \xi_{n+1} < \xi_{n+2}$, where the additional 6 knots are introduced in order to span the spline space ${\cS}_n$ of dimension $n$. 
\end{assumption}

We look for a cubic spline curve defined on the given dataset, 
whose definition incorporates both function and derivative constraints
to better design offsets and bi-offsets.

We assume to fit the data by a regression model partially based on  the same
assumptions that define a P-spline, a regression model introduced in \cite{EilersMarxFlexibleSmoothing}, and defined following \cite{CC2021}: 
\begin{definition} \label{def: unidim}
Given the {\bf Assumption} \ref{assump1}, the   P-spline  $s_{\lambda}(x)=\sum_{j=1}^{n}\alpha_j B_j(x)$ has coefficients that solve the penalized least square problem:
\begin{equation}\label{eq:Pspline_min}
 \min_{\alpha_1,\ldots, \alpha_{n}}\ 
\sum_{i=1}^{ {m}}  \left(y_i-\sum_{j=1}^{ {n}}  \alpha_j B_j(x_i)\right)^2\, +\,  {\lambda^2}  \,   {\sum_{j=3}^{n}\left((\Delta_2^{h} \bm{\alpha})_j \right)^2}
, 
\end{equation}
 with the  
{regularization parameter}  $\lambda>0$ balancing       data fidelity and regularization and  the  {penalty term}   defined by the difference operator on the coefficients   vector
$\bm{\alpha}$: 
\begin{equation}\label{eq:alphaPspline} 
(\Delta^{h}_2 \bm{\alpha})_{j} := \alpha_{j}  -2\alpha_{j-1}+\alpha_{j-2},\quad j=3,\ldots, n .  
\end{equation}
\end{definition} 
The  P-spline coefficients can also be seen as the solution of a Tikhonov regularization problem:
\begin{equation}\label{eq:Pspline}
\min_{\bm{\alpha}\in {\mathbb R}^{n}}\ 
 \| {\mathbf y} - {\mathbf B}\bm{\alpha}\|_2^2+\lambda^2  \|{\mathbf  D}\bm{\alpha}\|_2^2,     
\end{equation} 
where  $\mathbf{y}$ is the data vector of the observed data $\{y_i\}$, $i=1,\ldots,m$, ${\bf D}\in{\mathbb R}^{(n-2)\times n}$ is the difference  {three-banded} matrix:
  \begin{equation}
      {\mathbf D}{:}= \left[ \begin{array}{cccccc}
1       & -2  & 1        & 0       & \cdots      & 0 \\
  0        & 1       & -2  & 1     & \cdots      & 0 \\
\vdots     & \ddots             & \ddots     & \ddots  &  \cdots                  & \vdots\\
\vdots     & \vdots     & \vdots      & 1       & -2  & 1      \\
\end{array} \right],
\label{eq:1D_D}
\end{equation}
   the collocation matrix of the B-splines is 
    \begin{equation}
 \label{eq:1defB}
 {\bf B}\in {\mathbb R}^{m\times n}, \qquad 
 (\textbf{B})_{ij}=(B_i(x_{j})),
 \end{equation}  which depends on    $\{x_i\}$, $i=1,...,m$  and 
whose band structure is inherited by the locality of  B-splines.

  \subsubsection*{Step 1}  
  In this step,   we introduce a new regression spline defined as follows.

 \begin{definition}
    Given the {\bf Assumption} \ref{assump1}, the {\bf  Trajectory Penalized spline} (TP-spline)  $g(x)=\sum_{j=1}^{n}\alpha_j B_j(x)$ is a {\em modified   P-spline}, 
      whose coefficients  are computed by solving the following penalized least-square problem:
      \begin{equation}\label{eq:defTP}
       \min_{\bm{\alpha}\in {\mathbb R}^n}  \|{\mathbf y} - {\mathbf B}\bm{\alpha}\|_{2}^{2} 
        + \mu^2 \left(\|  {\mathbf \Phi}\bm{\alpha} - \Delta {\bf y}\|_{2}^{2} 
        + \|{\mathbf \Psi}\bm{\alpha} - {\bf z}\|_{2}^{2} \right)
        + \lambda^2  \|{\mathbf  D}\bm{\alpha}\|_2^2,     
    \end{equation}
  that is   \footnotesize{\begin{eqnarray}
        \min_{\alpha_1,\ldots ,\alpha_n} &&\sum_{j=1}^{N}  \left(\sum_{j=1}^{n}\alpha_j B_j(x_{j}) - y_{j}\right)^{2} 
        + \mu^2 \sum_{j=1}^{N} \left(\left(\sum_{j=1}^{n}\alpha_j B_j(x_{j})\right)' - \Delta y_{j}\right)^{2} 
        + \nonumber\\
        &+&\mu^2 \sum_{j=1}^{N} \left(\left(\sum_{j=1}^{n}\alpha_j B_j(x'_{j})\right)' - z_{j}\right)^{2} 
        + \,  {\lambda}^2  \,   {\sum_{j=3}^{n}\left((\Delta_2^{h} \bm{\alpha})_j \right)^2}\nonumber
    \end{eqnarray}}
    \normalsize
   with  $\bm{\alpha}$   the vector of the spline coefficients, $\lambda, \mu \in {\mathbb{R}}$   two regularization parameters,  $\textbf{D}$ as in (\ref{eq:1D_D}),  $\textbf{B}$ as in (\ref{eq:1defB}),  
    ${\bf \Phi}\in {\mathbb R}^{m\times n}$,
 with $({\bm \Phi})_{ij}=(B'_i(x_{j}))$.
      Moreover,   $\mathbf{x'}\in {\mathbb R}^{m-1}$  is the vector of the middle points  of $\mathbf{x}\in {\mathbb R}^{m}$,
and   $\mathbf{z} \in {\mathbb R}^{m-1}$  is  the     corresponding vector of the half differences of the $\mathbf{y} \in {\mathbb R}^{m}$ components;   finally
 ${\bf \Psi}\in {\mathbb R}^{(m-1)\times n}$,
with  $({\bm \Psi})_{ij}=(B'_i(x'_{j}))$.

\end{definition}

The meaning of $\mathbf{x'}$ and $\mathbf{z}$ will be made clear in the next section.    Moreover  $(\Delta {\bf y})_j =\frac{y_{j}-y_{j-1}}{x_{j}-x_{j-1}}$,  and $\Delta_2^{h} \bm{\alpha}$ is like in (\ref{eq:alphaPspline}).
 Differentiating  the (\ref{eq:defTP}) with  respect  to  ${\bm \alpha}$  
the least-square minimization problem 
  leads  to  the normal equations 
  \begin{align*}
\left( {\textbf{B}}^T{\textbf{B}} + \mu^2 \left( {\bm \Phi}^T{\bm \Phi} + {\bm \Psi}^T {\bm \Psi} \right) + \lambda^2 \textbf{D}^T \textbf{D} \right){\bm \alpha}= {\textbf{B}}^T {\textbf{y}} + \mu^2 \left( {\bm \Psi} {\textbf{z}} + {\bm \Phi}^T \Delta {\textbf{y}} \right).
\end{align*}
It is worth noting that the second order discretization is introduced to suppress in $g$   possible  small, sharp hooks which are dramatically enlarged in the offset definition. Such a hook is, in fact, related to large values in the second derivative. 
In this phase of the algorithm we do not use any quadrature methods. Instead, we approximate the integral term using the Eilers and Marx \cite{EilersMarxFlexibleSmoothing} idea, with the second-order finite difference operator defined  in (\ref{eq:alphaPspline}).

The regularization   parameters $\mu$ and $\lambda$ can be selected by a two-variable generalized cross-validation. This approach involves another minimization problem, in which the cost function is given by
\begin{equation}\label{eq:hGCV}
    h(\mu,\lambda)=\dfrac{\dfrac{1}{n}\Vert(I-\textbf{H}(\mu,\lambda))\textbf{y}\Vert_2^2}{\left(1 -\dfrac{1}{n}\text{tr}(\textbf{H}(\mu,\lambda))\right)^2}
\end{equation}
where 
\begin{equation*}
\textbf{H}(\mu,\lambda)=\textbf{B}\textbf{A}^{-1}\textbf{B}^T , 
\end{equation*}
 {with $\textbf{A}=  {\textbf{B}}^T{\textbf{B}} + \mu^2 \left( {\bm \Phi}^T{\bm \Phi} + {\bm \Psi}^T {\bm \Psi} \right) + \lambda^2 \textbf{D}^T \textbf{D}$}.
Once we have selected $\mu, \lambda$ as
\begin{equation}\label{gcv2d}
    (\mu^*,\lambda^*) = \mathop{\mathrm{argmin}}\limits_{\mu,\lambda>0} \quad h(\mu,\lambda)
\end{equation}
we can compute the vector ${\bm \alpha}$ by solving the linear system
obtained by minimizing the (\ref{eq:hGCV}).

     \subsubsection*{Step 2a}  
 Once the regression    model $g$ has been carried out, we   compute the two  offset curves.
 We cannot compute the offset ordinates directly according to the theoretical formulation. As already pointed out, the definition  leads to  cusps  and self-intersection points. To avoid these irregularities, we define a regression model defining the offset curve by
 solving the problem in the following definition:
 
 \begin{definition}
 \label{defOPspline}
  Given the {\bf Assumption} \ref{assump1},  the  {\bf Offset Penalized spline} (OP-spline)   $f$ approximating an offet curve, at a fixed distance $\tau>0$,
       from a function $g$, is a  regression  spline $f(x)=\sum_{j=1}^{n}\beta_j B_j(x)$   
  formulated by solving a   constrained optimization problem  
    \begin{align}\label{probStep2a}
  &\mathop{\mathrm{ min}} \mathcal{R}(\bm{\beta})  \\
&\,\,\text{s.t.}\quad f({\mathbf{x}}_o)=g(\overline{\mathbf{x}})+\tau g^\perp(\overline{\mathbf{x}})\nonumber \\
& \quad \quad \,\,\,\,g^\perp(\overline{\mathbf{x}})\cdot f'({\mathbf{x}}_o)=0 .\nonumber 
\end{align}
where 
  $\mathcal{R}(\bm{\beta})$  is  the functional defined by:
    \begin{equation*}
      \mathcal{R}(\bm{\beta}) = \int_{x_{o,\min}}^{x_{o,\max}} \Vert f''\Vert_2^2 \,\text{d}x, \qquad f''(x)=\sum_{j=1}^{n}\beta_j B''_j(x).
    \end{equation*}
    and the constraints are assumed to be at
    any set of     abscissae
$\bar{\mathbf{x}}=\{({\bar{{x}}}_j)\}$, $j=1,\dots,q$ 
with the corresponding points on the offsets along the point-wise  {orthogonal} directions,  given by
     \begin{equation*}
        \mathbf{x}_o^+ = \bar{\mathbf{x}} - \tau \dfrac{g'(\bar{\mathbf{x}})}{\sqrt{1 + g'(\bar{\mathbf{x}})^2}} .
    \end{equation*}
      \end{definition}

      It is worth noting that the methodology is analogous in both cases;  for simplicity of notation  we focus  on  the positive offset. For both, it is sufficient considering  
      
      \begin{equation*}
        \textbf{x}_o^\pm = \bar{\textbf{x}} \mp \tau \dfrac{g'(\bar{\textbf{x}})}{\sqrt{1 + g'(\bar{\textbf{x}})^2}}.
    \end{equation*}
    
As concern with the  two constraints in (\ref{probStep2a}):
    \begin{equation*}
        f(\textbf{x}_o) = g(\bar{\textbf{x}}) + \tau \dfrac{1}{\sqrt{1 + g'(\bar{\textbf{x}})^2}}
    \end{equation*}
     enforces that the offset curve lies at a signed distance $\tau$ from the generating curve, while
    \begin{equation*}
        g^{\perp}(\bar{\textbf{x}}) \cdot f'(\textbf{x}_o) = 0
    \end{equation*}
    is an equivalent formulation of the parallelism condition between $f^\perp(\textbf{x}_o)$ and $g^\perp(\bar{\textbf{x}})$.

     \subsubsection*{Step 2b (Refinement)}

 In order to deduce a functional  approximation of the original trajectory,
starting from information about the offsets, we focus on the assumption that, at any point of the offset curve, the unit tangent vectors  
 should  have consistent orientations. 
  We refer to this improvement  as "refinement" step, and  we aim to obtain an approximation of the first derivative of $f$ that satisfies two  conditions:
\begin{description}
    \item[(1)] The tangent vectors along the two curves maintain the same orientation;
    \item[(2)]  The tangent vector at each point of the offset becomes strongly proportional to the unit tangent vector at the corresponding point of the generator $g$.
\end{description}   

With this aim, let $\bar{\bar{\textbf{x}}}=\{(\bar{\bar{{x}}}_j)\}$, $j=1,\dots,p$ be a new set of abscissae, and let $\bar{\bar{\textbf{x}}}_o=\{{(\bar{\bar{{x}}}_0}_j)\}$, $j=1,\dots,p$  denote the     abscissae for the corresponding points on the offsets along the point-wise  {orthogonal} directions. Starting  from an initial  estimate  for $f$
given by Step 2a,  we search for a vector $\bm{\eta}$ such that:
\begin{equation*}
    f'(\bar{\bar{{x}}}_{o,j}) \approx \eta_j \dfrac{g'(\bar{\bar{{x}}}_j)}{\Vert g'(\bar{\bar{\textbf{x}}}) \Vert_2}, \quad j=1,\dots,p .
\end{equation*}
Given the coefficient vector $\bm{\beta}$ of $f$, we   compute the corresponding values of ${\bm \eta}$   by solving the following nonlinear least squares problem:
\\
 
    \begin{equation}\label{eq:pbStep3}
        \min_{\bm{\eta}} \sum_{j=1}^p \left(f'(\bar{\bar{{x}}}_{o,j}) - \eta_j \frac{g'(\bar{\bar{{x}}}_j)}{\|g'(\bar{\bar{\textbf{x}}})\|_2} \right)^2.
    \end{equation}
    
         \subsubsection*{Step 3}  
 
  In this final step, we aim to  approximate the generating TP-spline, given the OP-splines refined in Step 2b:
  \begin{definition} 
  Given the {\bf Assumption} \ref{assump1}, and let $f$ be the OP-spline of Definition \ref{defOPspline}.
  The {\bf {\em Bi-Offset}  spline} (BO-spline)   
  $h(x)=\sum_{j=1}^{n}\gamma_j B_j(x)$   is a {\em    regression spline}, computed  by
  solving the following nonlinear least squares problem:
\begin{equation}\label{eq:bioffset}
    \min_{\bm{\gamma}\in {\mathbb R}^n} \quad \Vert h(\bar{\bar{\mathbf{x}}}_o) - g(\bar{\bar{\mathbf{x}}})\Vert_2^2
\end{equation}
where
    $\bar{\bar{\mathbf{x}}}=\{(\bar{\bar{{x}}}_j)\}$, $j=1,\dots,p$ 
      is the set of points used to compute $\bm{\eta}$, and the  
            abscissae for the corresponding points on the bi-offset     $h$, along the point-wise  orthogonal directions   from the refinement, 
    are given by:
\begin{equation*}
\bar{\bar{\mathbf{x}}}_0=\bar{\bar{\mathbf{x}}} - \tau \dfrac{\eta ^T \cdot \dfrac{g'(\bar{\bar{\mathbf{x}}})}{\Vert g'(\bar{\bar{\mathbf{x}}})\Vert_2}}{\sqrt{1+\left( \eta ^T \cdot \dfrac{g'(\bar{\bar{\mathbf{x}}})}{\Vert g'(\bar{\bar{\mathbf{x}}})\Vert_2}\right)^2}}    .
\end{equation*}
 \end{definition}

\section{Algorithmic details and insights\label{sec:4}}
 
 In this section we describe some details  and computational insights in our approach.
 
As for the Step 1, in (\ref{eq:defTP})
a set of middle points $\mathbf{x}'$ is introduced    to adequately handle piecewise linear functions,  for  which  $\mathbf{x}$  is  chosen  as the  'kinks'  of  the  function (i.e. the points where the function is not smooth, meaning its derivative is not continuous) and 
\begin{equation*}
\mathbf{x}' = \frac{1}{2}
\begin{bmatrix}
1 & 1 &  \\
 & \ddots & \ddots & \\
& & 1 & 1 
\end{bmatrix}\mathbf{x}
\end{equation*}
At  anyone of   these  average values,  the  derivative  of  the  piecewise  linear  function  is  clearly  defined  as
\begin{equation*}
\mathbf{z} = \frac{1}{2}
\begin{bmatrix}
-1 & 1 &  \\
 & \ddots & \ddots & \\
& & -1 & 1 
\end{bmatrix}\mathbf{y}
\end{equation*}

This  is  simpler  and  more  stable  than  estimating  the  derivative  at  $\mathbf{x}$.
\\
 
As for the Step 2a,
the constrained minimization problem (\ref{probStep2a}) produces an improvement in terms of curvature and imposes a strong parallelism condition between the normal direction of the generator curve at $\bar{\textbf{x}}$ and the normal component of the offset curve at $\textbf{x}_o$. 
We solve (\ref{probStep2a}) writing  an equivalent form. For this purpose, it follows from the representation of $f$ that:
\begin{equation*}
	\Vert f'' \Vert_2^2 = \sum_{k,j} \beta_k^T \beta_j B^{''\, k T} B^{''\, j}.
\end{equation*}
Hence, the functional of the minimization problem is
\begin{equation*}
	\mathcal{R}(\bm{\beta}) = \int_{x_{\min}}^{x_{\max}} \Vert f''(\bar{\mathbf{x}})\Vert_2^2 \, d\bar{x} = \bm{\beta}^T \textbf{R} \bm{\beta}
\end{equation*}
where,  $\textbf{R}$ is a real symmetric matrix whose entries are given by
\begin{equation*}
\textbf{R}_k^{\;\; j} = \int_{x_{\min}}^{x_{\max}} B''(\bar{\mathbf{x}})^{k\,T} B''(\bar{\mathbf{x}})^j  d\bar{\mathbf{x}},
\end{equation*}
i.e., the integral of the product of the second derivatives of the corresponding basis functions.
As for the constraint $f(\overline{\textbf{x}}) = g(\overline{\textbf{x}}) + \tau g^\perp(\overline{\textbf{x}})$, 
the function $g^\perp$ represents the orthogonal component of the parametrized representation of the Cartesian curve $g(\overline{\textbf{x}})$.
Since the functions $g$ and $g^\perp$ evaluated at the set of points $\overline{x}$ are column vectors, the constraint can be written as
\begin{equation*}
\mathbf{B}(\bar{\mathbf{x}})\bm{\beta} = f(\bar{\mathbf{x}}) = g(\bar{\mathbf{x}}) + \tau g^\perp(\bar{\mathbf{x}}),
\end{equation*}
which is a linear system in the unknown $\bm{\beta}$.
\\
As for the constraint $g^\perp(\bar{\mathbf{x}}) \cdot f'(\bar{\mathbf{x}}) = 0$, 
using the matrix representations, this constraint can be expressed as:
\begin{equation*}
g^\perp(\bar{\mathbf{x}}) \cdot \mathbf{B}'(\bar{\mathbf{x}}){\bm \beta} = 0,
\end{equation*}
which represents a vanishing condition for the inner product of the normal vector of $g$ and the tangent vector of $f$. This is a linear equation in the unknown $\beta$.
Based on this, a compact formulation of  the constrained problem can be expressed as follows:
\begin{eqnarray}\label{eq:beta}
	&&\min_{{\bm \beta}} {\bm \beta} ^T \textbf{R} {\bm \beta}  \\
	&&\text{s.t.}\ \textbf{A}{\bm \beta}= \textbf{b}\nonumber
\end{eqnarray}
where
\begin{equation*}
A = \begin{bmatrix}
	\mathbf{B}(\bar{\mathbf{x}}) \\
	g^\perp(\bar{\mathbf{x}})^T \mathbf{B}'(\bar{\mathbf{x}})
\end{bmatrix}, \quad
b = \begin{bmatrix}
	g(\bar{\mathbf{x}}) + \tau g^\perp(\bar{\mathbf{x}}) \\
	0
\end{bmatrix}.
\end{equation*}
 \ \\
 
The refinement Step 2b  deserves extra attention. In the previous step, we imposed a strong parallelism condition between the tangent direction of the generator and the offset curves. However, to support a \textit{reversible} algorithm that allows us to define a derived  trajectory, starting from parallelism alone is not sufficient. Even though the two tangent components have the same direction, their orientation may change and their amplitudes are not proportional. In particular, a change in orientation introduces complications in the bi-offset reconstruction problem. The idea behind the refinement is to control these issues.

Starting from (\ref{eq:pbStep3}), set $T$   a diagonal matrix whose entries are the components of the unit tangent vector $\dfrac{g'(\bar{\bar{{x}}}_j)}{\Vert g'(\bar{\bar{\textbf{x}}}_j)\Vert_2}$, and given the matrix representation of the  spline $f$, the previous problem can be expressed in an equivalent form as
\begin{equation}
\label{min_norm}
    \min_{\bm\eta}\Vert \mathbf{B}'(\bar{\mathbf{x}}_{o}){\bm \beta}  - T{\bm \eta} \Vert_2^2,
\end{equation}
which is a quadratic minimization problem.

 To compute the vector ${\bm \eta}$, we adopt a Conjugate Gradient metohd.
As evidenced by the matrix formulation of the cost function:
\begin{equation*}
    \beta^T H_1 \beta - 2 \beta^T f \eta + \eta^T H_2 \eta
\end{equation*}
with $\mathbf{H}_1 = \mathbf{B}'^{\,T}\mathbf{B}'$, $\mathbf{H}_2 = \mathbf{T}^T \mathbf{T}$, and $\mathbf{f} = \mathbf{B}'^{\,T}\mathbf{T}$,
the optimization problem is quadratic and convex, with $H_2$ symmetric and definite positive matrix.
Consequently, the minimization of this function is equivalent to finding the stationary point where the gradient vanishes. 
A fundamental property of this method is that, by searching along a sequence of conjugate directions, it is guaranteed to converge to the exact solution in a finite number of iterations, while it is, of course, well-known that in many cases a truely iterative application of the method often leads to better results.

In the final step, we   carry  out the reconstruction of the original generating curve from its offset. The main idea of this phase is to construct a new offset curve in which the upper/lower
 offset becomes the new generator curve. In particular, for the purpose of returning to the original curve $g$ from the upper offset $f$, we need to impose the computation of its under offset using the same fixed non-zero radius $\tau$. 
 As   previously remarked, it is not possible to proceed using the classical theoretical approach, thus, we     engender a new regression minimum problem.

To this purpose we look for a new set of interpolating points by benefiting from the previous samples used in the refinement stage. Once computed     ${\bm \eta}$   from (\ref{min_norm})
it follows that the new set of interpolating points is given by:
\begin{equation*}
    \bar{\bar{x}}_{o,j} = \bar{\bar{x}}_j - \tau \dfrac{ \dfrac{\eta_j  g'(\bar{\bar{x}}_j)}{\Vert g'(\bar{\bar{x}}_j)\Vert_2} }{\sqrt{1+\left( \dfrac{\eta_j  g'(\bar{\bar{x}}_j)}{\Vert g'(\bar{\bar{x}}_j)\Vert_2} \right)^2}} ,
\end{equation*}
and the Step 3  is reduced to solving   the   non-linear least-squares problem:
\begin{equation*}
{{\bm \gamma}}^*= arg\min_{{\bm \gamma}} \Vert h(\bar{\bar{\mathbf{x}}}_o) - g(\bar{\bar{\mathbf{x}}}) \Vert_2^2=arg\min_{{\bm \gamma}}\Vert \bar{\bar{\mathbf{B}}}{\bm \gamma} - g(\bar{\bar{\mathbf{x}}})\Vert_2^2.
\end{equation*}
with  $\bar{\bar{\mathbf{B}}}:=B(\bar{\bar{\mathbf{x}}}_o)$   the collocation  matrix for  $h$.  
It is not essential to convert this unconstrained minimum problem into a new form, it can be solved by   the Matlab function \texttt{lsqnonlin}.
 The key element, which must be emphasized, is that the purpose of the refinement is concretely understood in the construction of the vector ${\bm \eta}$. This vector plays a fundamental role in the reconstruction, as it enables the construction of the new sampling - and thus the effective reconstruction - without relying on the derivatives of the  released spline $f$, but instead adopting those of the function $g$. This modification is crucial, since by enforcing the derivative of $g$ in the actual computation we prevent possible orientation inversions of $f$, while still preserving the same amplitude due to the carefully obtained proportionality factor.  
 {Algorithm \ref{algOPsplineTOTAL} synthesizes the whole 4-step procedure.}
 
 \begin{algorithm}[!h]
 \caption{TP-spline definition  and bi-offset reconstruction}
 \label{algOPsplineTOTAL} 
  \begin{algorithmic}[1]
\State\underline{Input}:
\State ${\mathcal D}=\{(x_i,y_i)\}$, $i=1,\ldots ,m$ set of     (noisy) data with $x_i \in [a,b]$
\State $\Xi=\{\xi_{-1},\ldots, \xi_{n+2}\}$, with  $\xi_{i}<\xi_{i+1}$, uniform (augmented) set of  knots  with $\xi_{2}\equiv a=x_1$, $\xi_{n-1}\equiv b=x_m$.
\State $B_1(x),B_2(x),\hdots, B_n(x)$ {cubic B-splines} defined on $\Xi$ 
 \State  \hspace{2mm} Step 1: Compute the TP-spline $g(\mathbf{x})$,  by   (\ref{eq:defTP}), with $(\mu^*,\lambda^*)$ by  (\ref{gcv2d}).
\State  \hspace{2mm} Step 2a: Construct the offsets  $f$ by (\ref{probStep2a})
\State  \hspace{2mm} Step 2b: Refine the offsets tangent vectors by (\ref{eq:pbStep3}):
\State  \hspace{2mm} Step 3: Construct  the bi-offset function by (\ref{eq:bioffset}).  
\State \underline{Output}:    $h= \sum_{j=1}^{n}\gamma_j B_j$ 
be the bi-offset curve of $g$.
\end{algorithmic}
   \end{algorithm}

The main issue encountered in the implementation of this algorithm lies in the non-unique and non-adaptive process for the selection and distribution of knots. Utilizing the sample data for the regression processes, the number of equidistant knots was set within a range of $20\% - 40\%$ of the sample data. This approach often produced trends that were similar to the model but exhibited undesired characteristics. We emphasize that the problem does not seem inherent to the model itself, as in regions where the knots are adequate, the process is correctly executed.

An adaptive selection of knots would help make the process more streamlined, functional, and easily adaptable from one curve to another, without the need to manually set the knots correctly through a process of trial and error.
Nevertheless, the penalty terms in the regression models, are generally useful for relaxing the importance of the location of the knots and their number.

\section{Numerical Experiments\label{sec:5}}



This section is to show our approach in action on different datasets.  
The experiments want to highlight the model reliability in the bi-offset design, the improvement given by the refinement strategy, and a comparison with classical models, available in literature.

{ 
The
tests were carried out on an 13th Gen Intel(R) Core(TM) i7-13620H (2.40
GHz), 3500 MHz, 2 cores, 4 logical processors. The MATLAB$^{\copyright}$ release is R2025b.}

The tests reported below are performed by varying the offset distance ${ \tau}\in 
\{0.1, 0.3, 0.5, 0.7\}$.  For each radius,  the theoretical
offsets, the offsets constructed according to our model and  the reconstruction of the base curve, with and without the application
of the refinement, are referred to.

We assume to have a dataset consisting of function values and first derivatives, associated to a sequence of given abscissas.
 In order to simulate noisy data, zero-mean Gaussian noise is added to both the function values and the derivatives with standard deviation $\sigma$ suitable specified.
The regularization parameters $(\mu, \lambda)$  are computed by solving a two-parameter minimization problem, where the objective function is in (\ref{eq:hGCV}). 
The (\ref{gcv2d}) is solved using the Matlab \texttt{fmincon} function as a local optimizer.
The minimization problem in (\ref{min_norm}) is solved using the  Matlab
function \texttt{pcg}.



\subsection{Test 1 (model reliability)}
\label{test:1}
In this section the experiments are to prove the model reliability,
by computing the accuracy in the offset definition and in the reconstruction of the base function through  the bi-offset.
We consider a dataset $(x_i,y_i)$, $i=1,\ldots,47$, with uniformly distributed  abscissae and $y_i=p_1(x_i)$ with
$$p_1(x)=\vert \sin (x) \cos(2x)  \vert,\, x\in [0,2\pi].$$
 We define the TP-spline using 
$n=14$ knots,  and   the computed regularization parameters $\mu=2.4628\times 10^{-2}$ and $\lambda =2.0506\times 10^{-2}$.
 
\begin{figure}[h!]
    \centering
    \includegraphics[width=\linewidth]{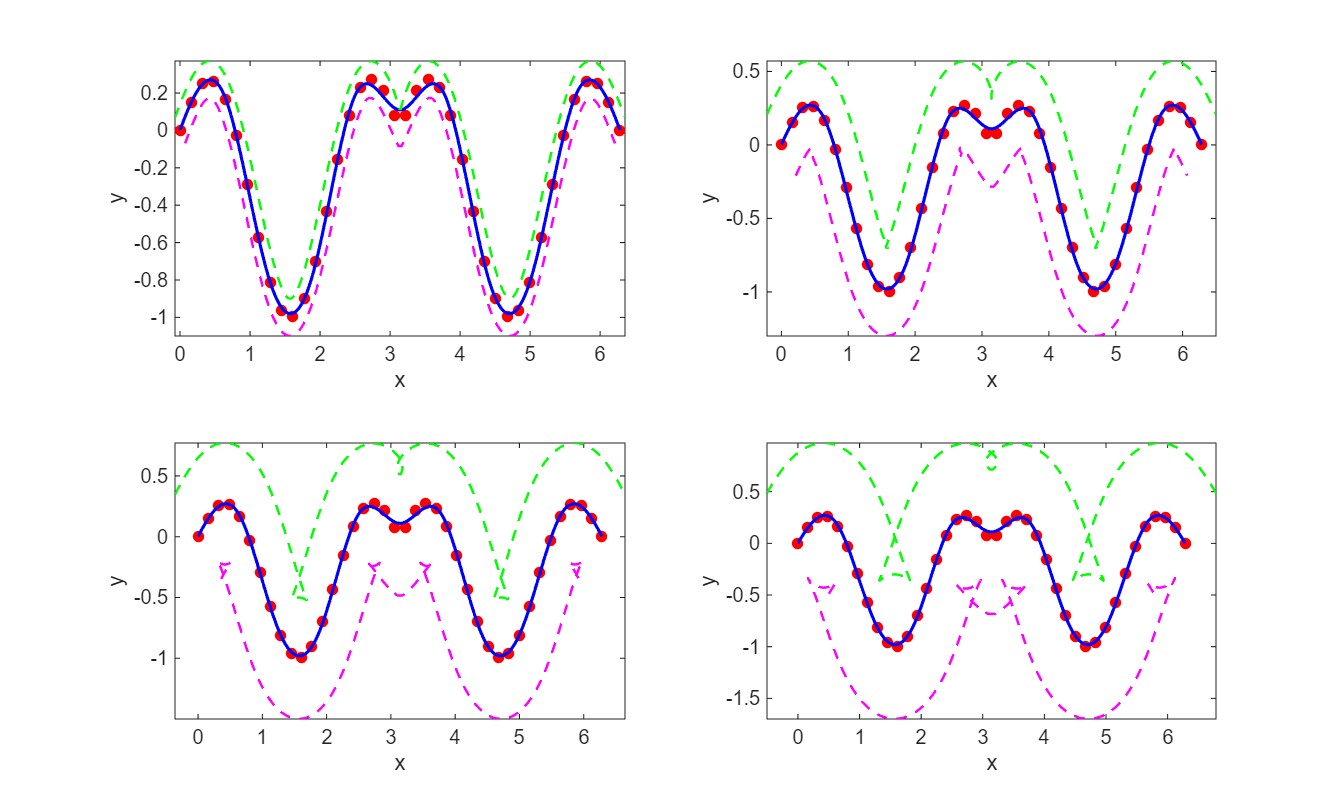}
    \caption{
     Base  curve  (blue '$\textcolor{red}{\bullet} -$') and theoretical offsets, {interior offset}   (magenta '- -') and  {exterior offset}   (green '- -'), for different ${ \tau}$: top left $\tau=0.1$,   right $\tau=0.3$;  bottom left $\tau=0.5$,    right $\tau=0.7$.\label{figTest1_1}}
\end{figure}

In Figure \ref{figTest1_1}  we describe  the theoretical offsets   according to the Definition \ref{def:offsets}, for different distances $\tau$. These results  highlight  the presence of cusps  and self-intersections.\\

\begin{figure}[h!]
    \centering
    \includegraphics[width=\linewidth]{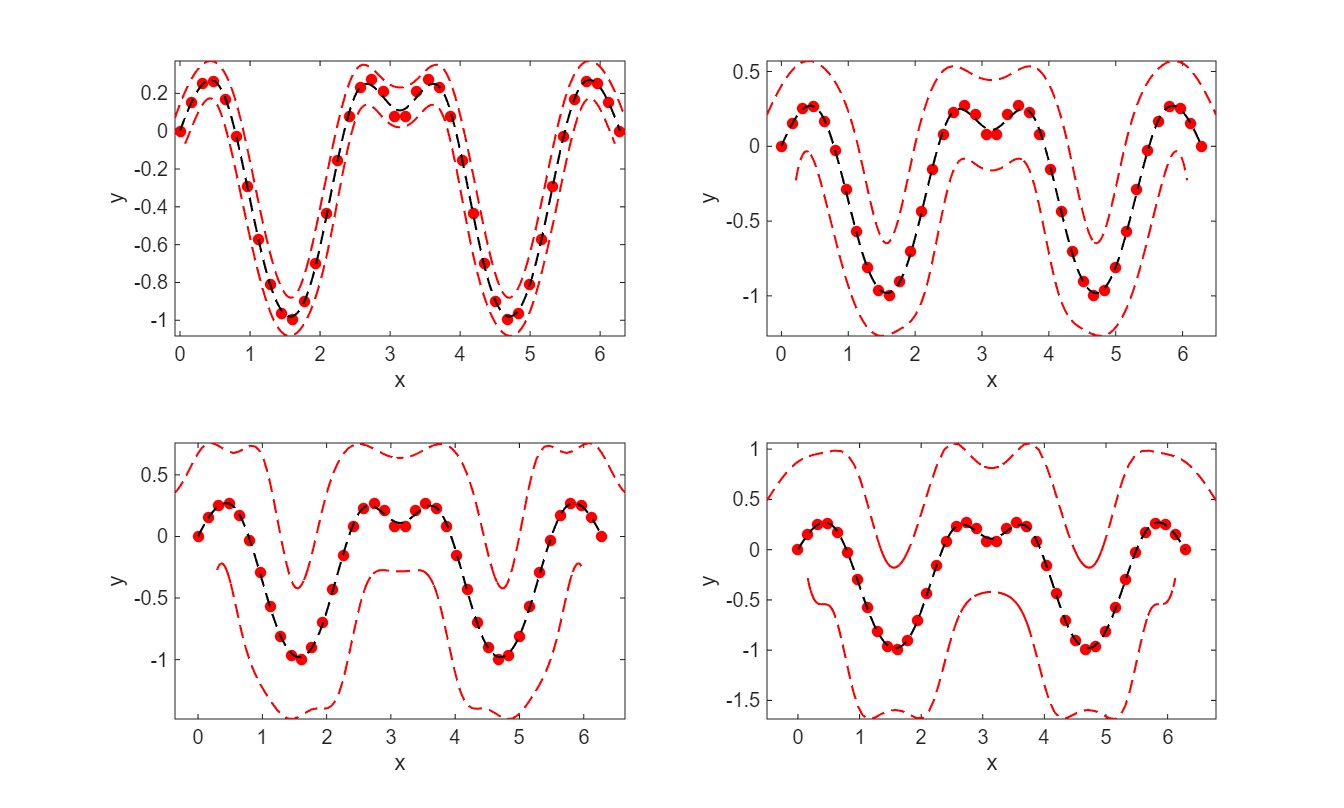}
    \caption{
     Base  curve  (black '$\textcolor{red}{\bullet} -$') and OP-splines  (red   '- -') for different ${ \tau}$: top left $\tau=0.1$,   right $\tau=0.3$;  bottom left $\tau=0.5$,    right $\tau=0.7$. \label{figTest1_2}}
 \end{figure}

\begin{figure}[!bth]
  \centering
		\subfigure[]{\includegraphics[width=0.8\textwidth]{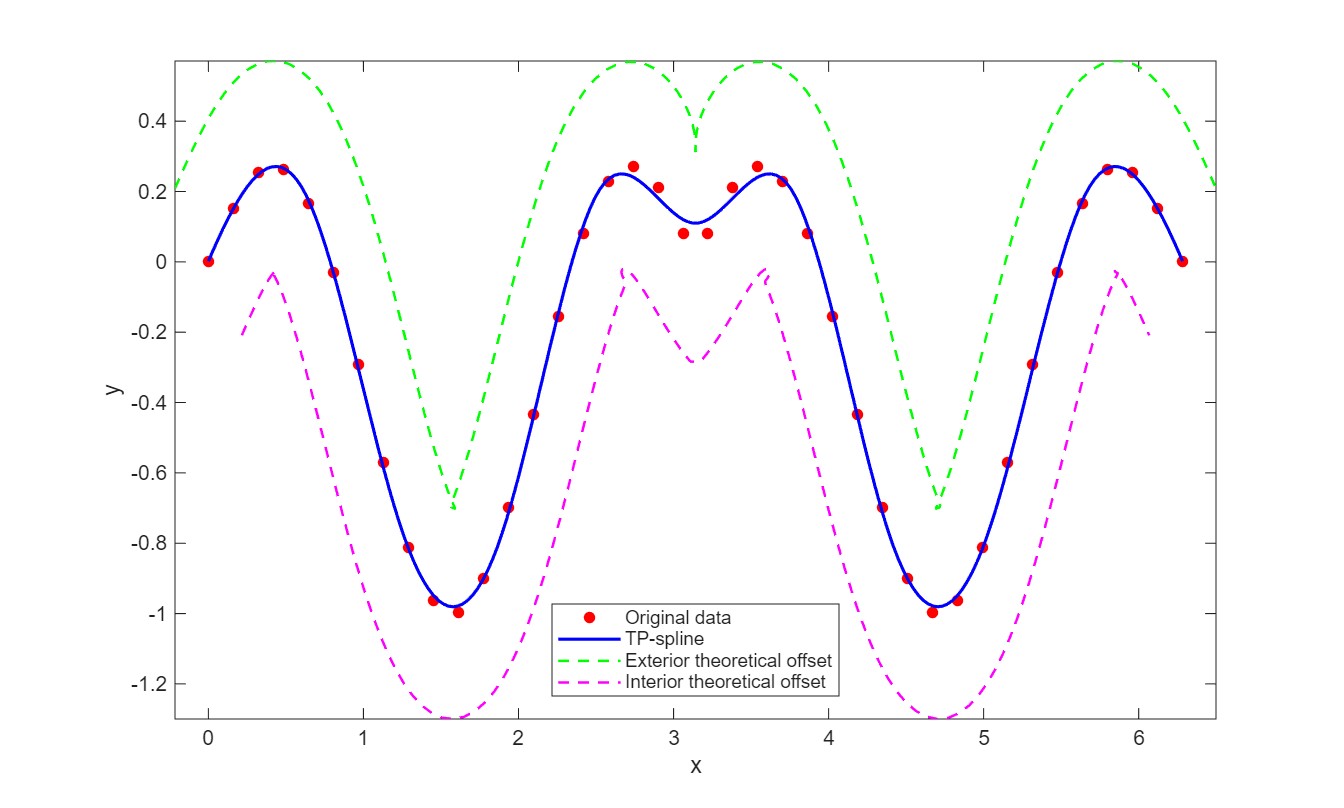}  }\\
		\subfigure[]{\includegraphics[width=0.8\textwidth]{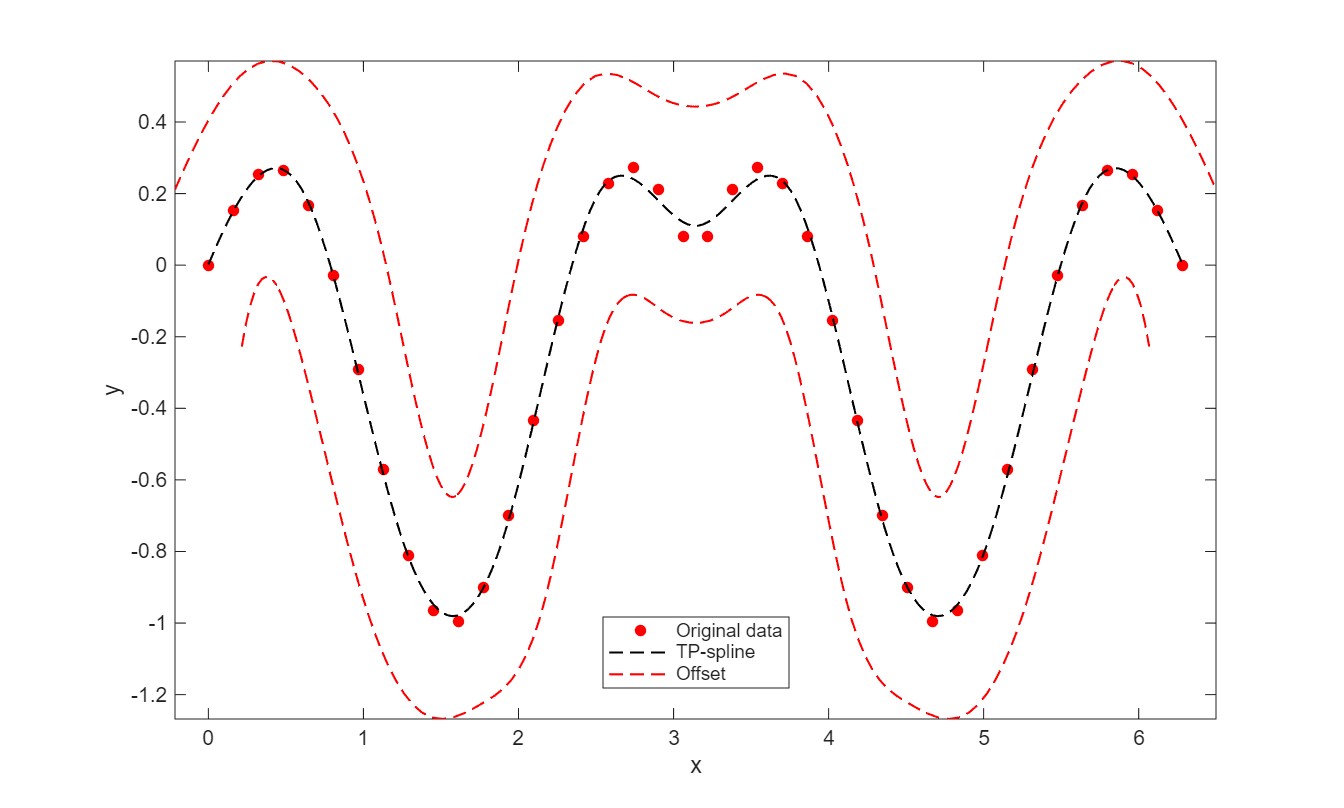}}
  \caption{Case $\tau=0.3$:    {interior}   (magenta '- -') and  {exterior}   (green '- -')  theoretical offsets (a),   OP-splines (red '- -')  (b).}
\label{figTest1_3}
\end{figure}

 In Figure \ref{figTest1_2}  the two OP-splines for different distances $\tau$ are presented.
In Figure 
\ref{figTest1_3}  we refer the comparisons between the two aforementioned approaches,   for a   specific  distance $\tau=0.3$.

\begin{figure}[!bth]
  \centering
		\subfigure[]{\includegraphics[width=0.8\textwidth]{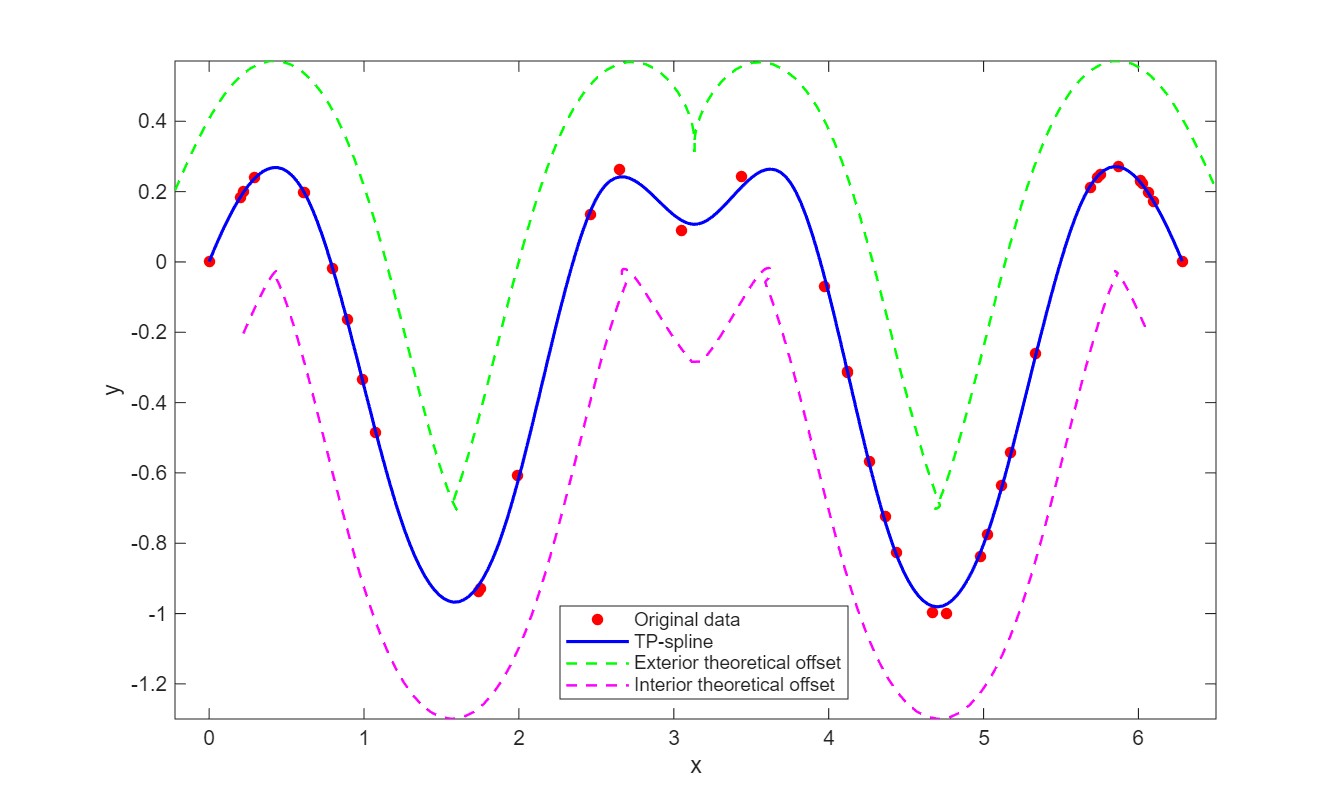}    }\\
		\subfigure[]{\includegraphics[width=0.8\textwidth]{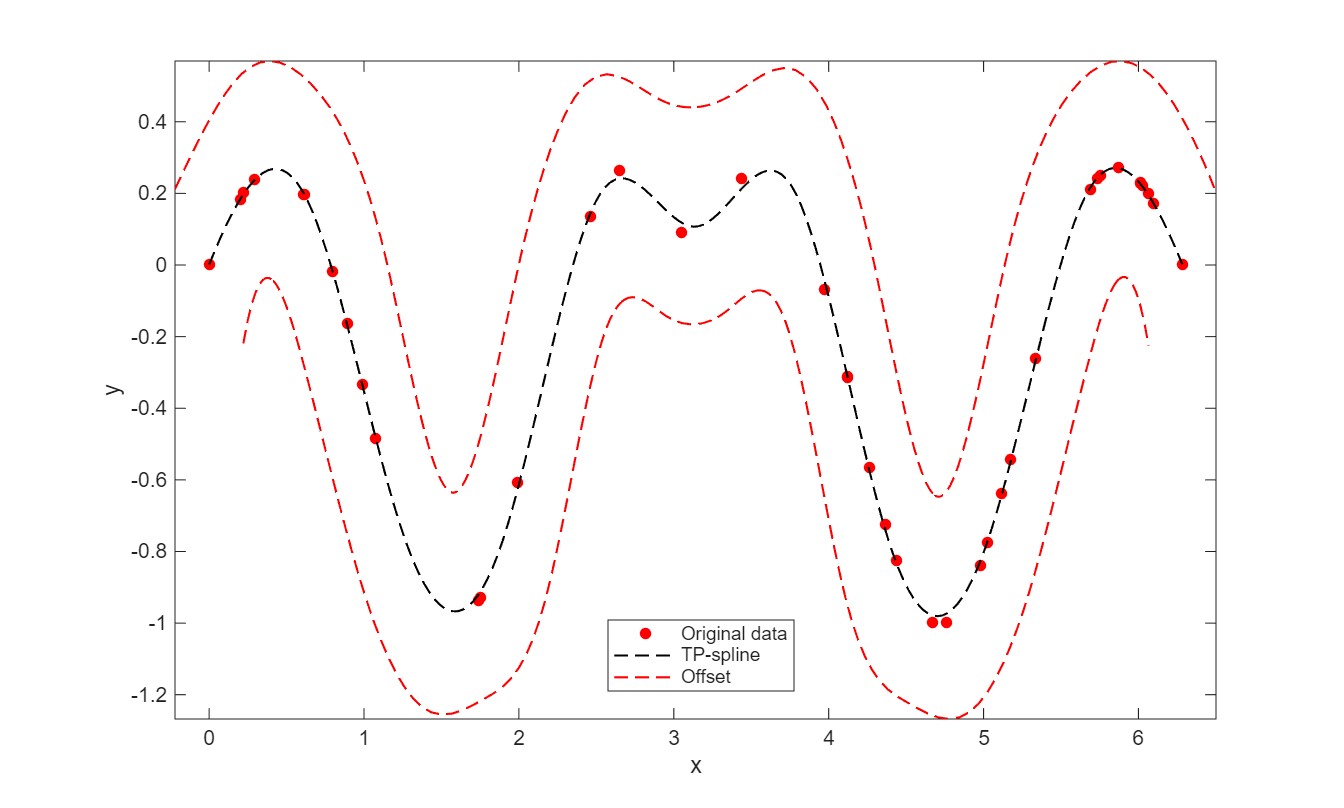}}
  \caption{Case $\tau=0.3$, $\sigma=1.0e-2$:    {interior}   (magenta '- -') and  {exterior}   (green '- -')  theoretical offsets (a),   OP-splines (red '- -')  (b).}
 \label{figTest1_4}
\end{figure}

 Similar results are obtained assuming   data affected by relative zero-mean Gaussian noise with standard deviation $\sigma=10^{-2}$. In this case  the computed regularization parameters  of the TP-spline are $\lambda= 1.1086 \times 10^{-2}$ and 
$\mu=2.9283\times 10^{-2}$ (see Figure \ref{figTest1_4}).

\subsection{Test 2 ({Refinement} or \emph{not} refinement)}
\label{test:2}
 
\begin{figure}[!bth]
  \centering
		\subfigure[]{ \includegraphics[width=0.8\linewidth]{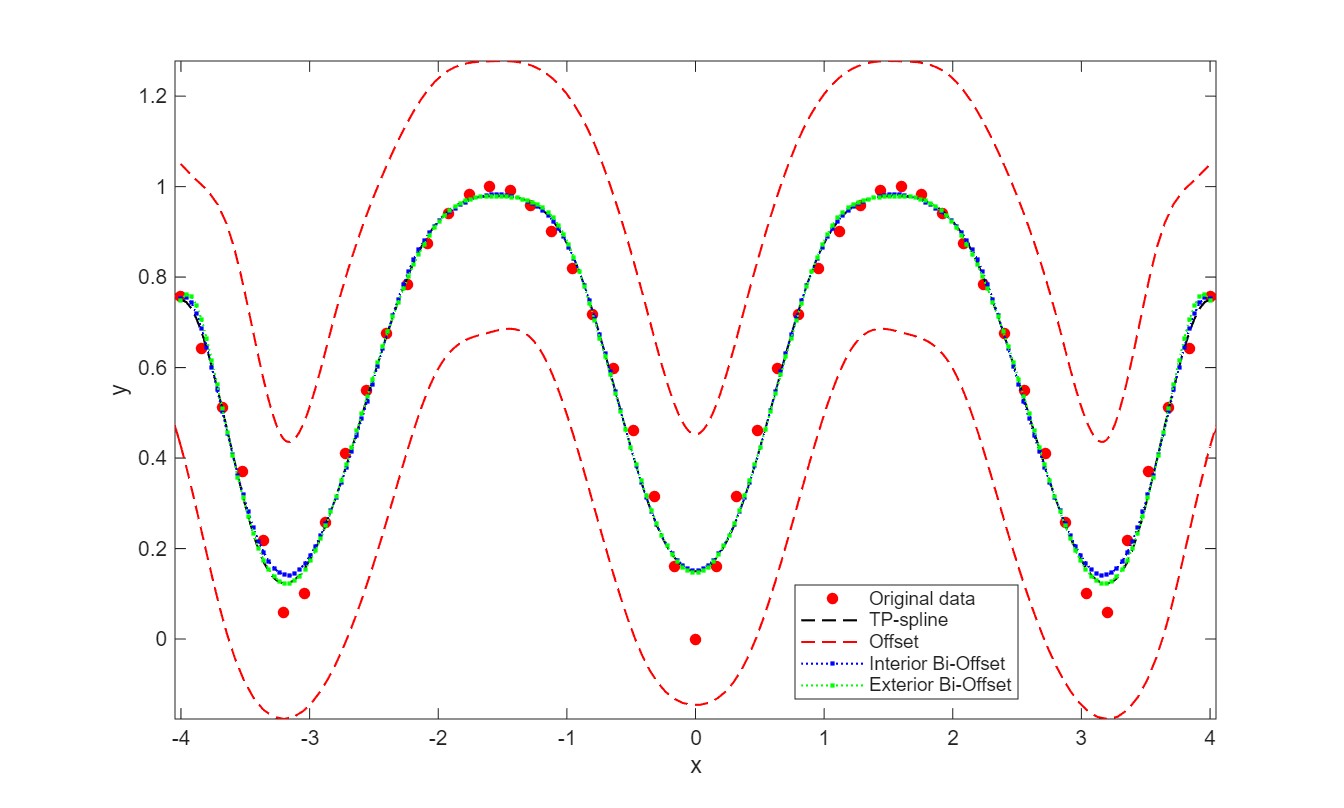}    }\\
		\subfigure[]{\includegraphics[width=0.8\linewidth]{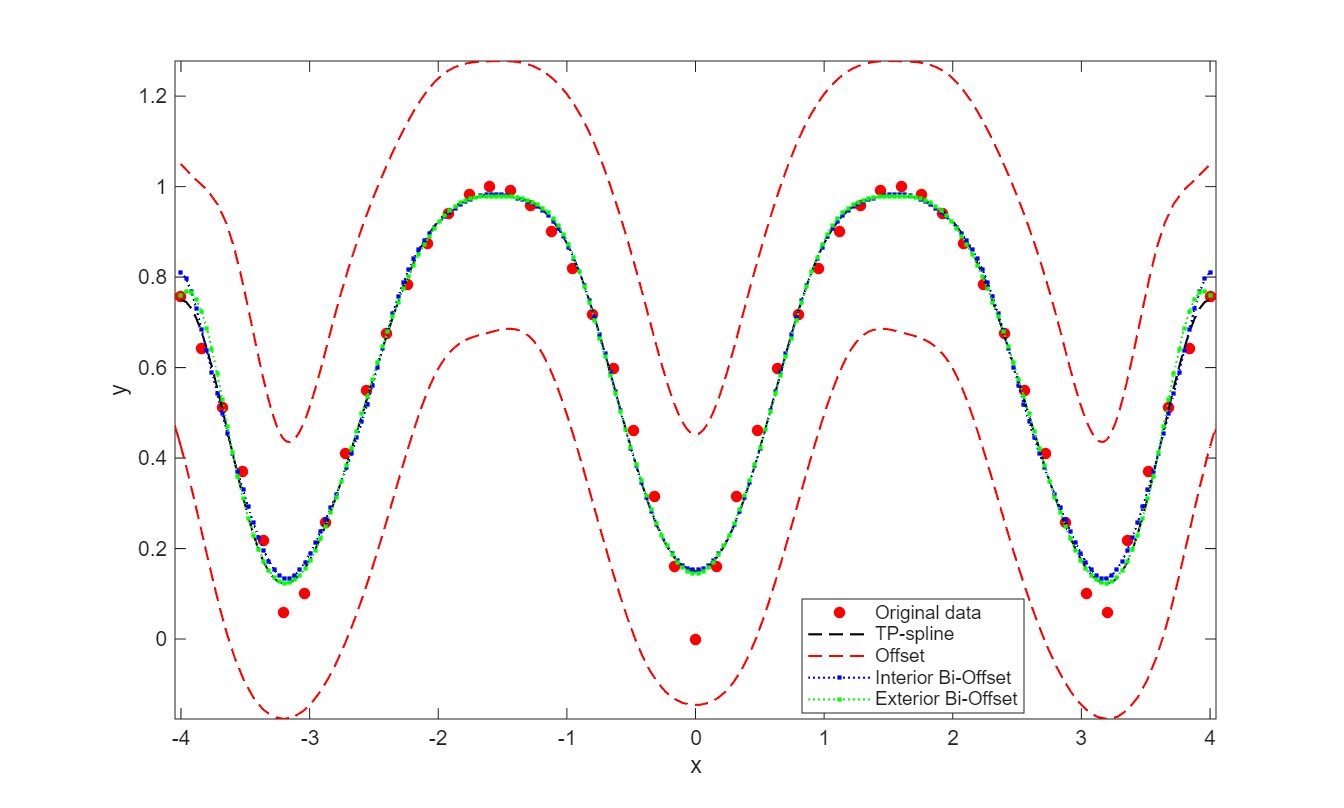}}
  \caption{Case $\tau=0.3$:   {interior (or from above)}   (blue '- -') and  {exterior (or from below)}   (green '- -')  BP-spline  with refinement   (a);    {interior (or from above)}   (blue '- -') and  {exterior (or from below)} BP-spline  without refinement (b).}
\label{figTest2_3}
\end{figure}

 We implement this test with the aim of highlighting the improvement in bi-offsets construction using the refinement.  We consider a dataset 
 $(x_i,y_i)$, $i=1,\ldots,51$, with 51  uniformly distributed  abscissae, 
  $y_i=p_2(x_i)$ with
 $$p_2(x)=\vert \sin (x)\vert,\, x\in [-4,4].$$ We define the TP-spline using
$n=14$ knots, and the computed regularization parameters  $\mu=4.3242\times 10^{-1}$ and $\lambda =3.6628\times 10^{-3}$ and compute the two BP-splines 
representing the two bi-offsets, respectively from below (\textit{exterior bi-offset}) and from above (\textit{interior bi-offset}).  Figure \ref{figTest2_3}
shows that the  refinement allows to manage the boundary effects arising in the bi-offset approximations.
Numerical results confirming these improvements are detailed in the tables of the  next Test section.

  

 \subsection{Test 3 (model comparison)}
This section presents some results on a comparison  between the bi-offset approximations
obtained both using our approach and     different spline models,  in the approximation of the dataset (in Step 1). Particularly we compare the  results obtained by  TP-spline, with the ones given by the classical    P-spline \cite{EilersMarxFlexibleSmoothing}, computed by   a proprietary code, and  the  interpolating cubic spline, by the MATLAB function \texttt{csaps}.\\
In the same contexts of Test \ref{test:1} and Test \ref{test:2}, we compute  
the \textit{Mean Squared Error} (MSE),
    quantifying the average squared deviation, 
    between  either  the interior or the exterior  bi-offsets with the aforementioned 
    approximating splines,
at different distances $\tau=\{0.1,0.3,0.5,0.7\}$.
The corresponding results are in the Tables \ref{Tab2}-\ref{Tab3} for the generating functions $p_1(x)$ and $p_2(x)$ respectively.  
We can observe that our model 
produces a better accuracy with respect to the compared spline  models,
when  the approximation is given  both by interior and by exterior bi-offset.
Furthermore these results improve when the refinement is considered.

{\small \begin{table}[h!]
\centering
 \begin{tabular}{|c|cc|cc|c|}
\hline
\hline
model &   Interior   & Interior   &   Exterior   & Exterior   & ${\tau}$ \\
  & \textit{with refinement}     &     & \textit{with refinement}      &     &   \\
\hline
{'TP-spline'} & 5.9752e-06 & 5.9072e-06 & 6.7020e-06 & 6.7080e-06 & 0.1 \\ 
{'P-spline' } & 5.9693e-06 & 5.9021e-06 & 6.6965e-06 & 6.7024e-06 & 0.1 \\
{'Spline'   } & 1.2617e-04 & 1.2604e-04 & 1.2112e-04 & 1.2126e-04 & 0.1 \\
\hline
{'TP-spline'} & 4.1069e-04 & 4.2294e-04 & 5.7735e-04 & 6.7511e-04 & 0.3 \\
{'P-spline' } & 4.1045e-04 & 4.2268e-04 & 5.7696e-04 & 6.7481e-04 & 0.3\\
{'Spline'   } & 6.1803e-04 & 6.0767e-04 & 1.0242e-03 & 1.2318e-03 & 0.3 \\
\hline
{'TP-spline'} & 1.8491e-02 & 2.6425e-02 & 9.6892e-03 & 1.5695e-02 & 0.5 \\
{'P-spline' } & 1.8485e-02 & 2.6382e-02 & 9.6921e-03 & 1.5707e-02 & 0.5 \\
{'Spline'   } & 2.7309e-02 & 3.0080e-02 & 1.3924e-02 & 2.6333e-02 & 0.5\\
\hline
{'TP-spline'} & 1.3061e-01 & 1.0827e-01 & 4.6631e-01 & 5.1531e-01 & 0.7\\
{'P-spline' } & 1.3055e-01 & 1.0822e-01 & 4.6798e-01 & 5.1790e-01 & 0.7 \\
{'Spline'   } & 1.2021e-01 & 1.2386e-01 & 2.6408e-01 & 9.2787e-01 & 0.7 \\
\hline
\hline
\end{tabular}
 \caption{Test function $p_1(x)=\vert \sin (x) \cos(2x) \vert$: Mean Square Error between   the specified spline models and the computed interior and exterior bi-offsets, with and without refinement, for different $\tau$ values. \label{Tab2}}
\end{table}
}

{\small 
 \begin{table}[h!]
\centering
\begin{tabular}{|c|cc|cc|c|}
\hline
\hline
model &   Interior   & Interior   &   Exterior   & Exterior   & ${\tau}$ \\
  & \textit{with refinement}     &     & \textit{with refinement}      &     &   \\
\hline
 
{'TP-spline'} & 1.0992e-06 & 1.4281e-06 & 7.0177e-07 & 6.3640e-07 & 0.1 \\
 
{'P-spline' } & 1.0721e-06 & 1.3217e-06 & 6.9827e-07 & 6.2982e-07 & 0.1 \\
 
{'Spline'   } & 1.6819e-04 & 1.6824e-04 & 1.5222e-04 & 1.5228e-04 & 0.1 \\
 
\hline
 
{'TP-spline'} & 5.1941e-05 & 1.5621e-04 & 4.4498e-05 & 1.1864e-04 & 0.3 \\
 
{'P-spline' } & 7.5100e-05 & 2.5420e-04 & 3.1850e-05 & 1.0551e-04 & 0.3 \\
 
{'Spline'   } & 2.3808e-03 & 2.3764e-03 & 1.5170e-04 & 1.5198e-04 & 0.3 \\
 
\hline
 
{'TP-spline'} & 5.7482e-03 & 1.7173e-01 & 3.3926e-04 & 8.2687e-04 & 0.5\\
 
{'P-spline' } & 2.2517e-01 & 8.9925e-02 & 4.0532e-04 & 8.8120e-04 & 0.5\\
 
{'Spline'   } & 1.7940e-02 & 1.7965e-02 & 1.8115e-04 & 1.8096e-04 & 0.5\\
 
\hline
 
{'TP-spline'} & 2.3323e-02 & 3.1002e-02 & 9.0385e-04 & 2.4493e-03 & 0.7 \\
 
{'P-spline' } & 3.8280e-02 & 3.0820e-02 & 1.1068e-03 & 2.0685e-03 & 0.7 \\
{'Spline'   } & 6.9913e-02 & 6.9978e-02 & 3.6229e-04 & 3.6269e-04 & 0.7  \\
\hline
\hline
\end{tabular}
\caption{Test function $p_2(x)=\vert \sin (x) \vert$: Mean Square Error between   the specified spline models and the computed interior and exterior bi-offsets, with and without refinement, for different $\tau$ values. \label{Tab3}}
\end{table}
}

\section{Conclusions}
In this paper we present a regularization technique   which incorporates both function values and 
derivatives, for designing   offset curves  and reconstruct an approximation of the base trajectory, given its offsets
  (we say the bi-offsets, for short).
An element of interest in the work carried out and 
worth to be dealt with
so far is undoubtedly our model's ability to reconstruct the original curve from its offsets. The comparison metrics introduced, using the mean square error, serve to confirm the good reliability of the developed process, attesting to its functionality.
Further relevance must be given to the results concerning  the relationship between the offset distance  $\tau$ and the  cusps presence,  connected to the curvature-radius relationship. 
 We verify empirically that for high values of the radius, the presence of cusps increases, so attention is due to   the correct choice of ${\tau}$, and
drives our interest  
 for future investigations, particularly because the selection of an {{optimal confidence interval}} or bounds  for the radius, 
concern  not only the {construction of the offsets} but also a regular and faithful possible {extraction of   road's external
 boundaries}  and {adaptive reconstruction of the   center line}, arising in the applications
 of automatic or semi-automatic assisted driving \cite{HU2018482,9703250}.

 \section*{Acknowledgements}
RC and SM are members of the INdAM research group GNCS, which has
partially supported this work. This research has been accomplished within RITA (Research ITalian network
on Approximation) and UMI-TAA groups. RC was partially supported by the Italian MUR through the PRIN2022 Project 'Numerical Optimization with Adaptive Accuracy and Applications to Machine Learning',   code: 2022N3ZNAX,  and 
the PRIN2022-PNRR Project 'A multidisciplinary approach to evaluate ecosystems resilience under climate change',   code: P2022WC2ZZ.

 \section*{Declarations}
The MATLAB code used in this study is available from the authors upon request.

 \end{document}